 \documentclass[10pt]{article}
 \usepackage{amsmath}
 \usepackage{pifont}
 \usepackage{amsfonts}
 \usepackage{mathrsfs}

 \usepackage{booktabs}
 \usepackage{threeparttable}
 \usepackage{mathrsfs,amsfonts,amsmath}
 \usepackage{harvard}
 \usepackage[dvips]{color}

 \newtheorem{theorem}{Theorem}[section]
 \newtheorem{definition}[theorem]{Definition}
 \newtheorem{lemma}[theorem]{Lemma}
 \newtheorem{corollary}[theorem]{Corollary}
 \newtheorem{proposition}[theorem]{Proposition}
 \newtheorem{remark}[theorem]{Remark}
 \newtheorem{condition}[theorem]{Condition}
 \newtheorem{example}{Example}[section]

 \def\blemma{\begin{lemma}}\def\elemma{\end{lemma}}
 \def\bproposition{\begin{proposition}}\def\eproposition{\end{proposition}}
 \def\ttheorem{\begin{theorem}}\def\etheorem{\end{theorem}}
 \def\bcorollary{\begin{corollary}}\def\ecorollary{\end{corollary}}
 \def\bremark{\begin{remark}}\def\eremark{\end{remark}}
 \def\bcondition{\begin{condition}}\def\econdition{\end{condition}}

 \def\benumerate{\begin{enumerate}}\def\eenumerate{\end{enumerate}}
 \def\bitemize{\begin{itemize}}\def\eitemize{\end{itemize}}

 \def\beqlb{\begin{eqnarray}}\def\eeqlb{\end{eqnarray}}
 \def\beqnn{\begin{eqnarray*}}\def\eeqnn{\end{eqnarray*}}
 \def\ar{\!\!\!&}

 \def\mbb{\mathbb}\def\mbf{\mathbf}

 \def\qed{\hfill$\Box$\medskip}

\begin{document}

\noindent{(Version: 2013/04/11)}

\bigskip\bigskip

\centerline{\Large\bf Parameter Estimation in Two-type Continuous-}

\smallskip

\centerline{\Large\bf state Branching Processes with Immigration}

\bigskip

\centerline{Wei Xu}

\medskip

\centerline{School of Mathematical Sciences, Beijing Normal University}

\medskip

\centerline{Beijing 100875, People's Republic of China}

\medskip

\centerline{E-mails: \texttt{xuwei@mail.bnu.edu.cn}\quad\textbf{\textit{or}}\quad\texttt{xuweijlu.2008@163.com}}
\bigskip

{\narrower

\noindent\textit{Abstract:} We study the estimation of two-type
continuous-state branching processes with immigration (CBI-processes).
The ergodicity of the processes is proved. We also establish the strong
consistency and central limit theorems of the conditional least squares
estimators and the weighted conditional least squares estimators of the
drift and diffusion coefficients based on low frequency observations.

\smallskip


\par}


\section{Introduction}

\setcounter{equation}{0}

Branching processes have been used widely not only in biology, but also in
financial world. For example, Galton-Watson branching processes with
immigration (GWI-processes) are used to study the evolution of different
species. Continuous-state branching processes (CB-processes) were first
introduced by  Ji\v{r}ina (1958). In particular, a continuous CB-process
can be obtained as the unique solution of a stochastic equation system
driven by Brownian motion. Kawazu and Watanabe (1971) constructed
continuous-state branching processes with immigration (CBI-processes). In
view of the results of Dawson and Li (2006), a general single-type
CBI-process is the unique strong solution of a stochastic equation driven
by Brownian motions and Poisson random measures. The two-type CB-processes
was first be introduced by Watanabe (1969). Ma (2012) proved the existence
and uniqueness of the strong solution of a two-dimensional stochastic
integral equation system with jumps. He also showed that the unique
solution is a two-type CBI-process. In financial world, multitype
CBI-processes are used to describe the relations of the prices of
different assets and interest rates of different currencies.

Firstly, we introduce a special continuous single-type CBI-process defined
by the following equation:
 \beqlb\label{1.1}
dX_t=(a-bX_t)dt+\sigma\sqrt{X_t}dB_t,
 \eeqlb
where $(a, b, \sigma)\in (0, +\infty)^3$ and $B_t$ is a standard Brownian
motion. In fact, the solution of (\ref{1.1}) is also known as the
\textit{Cox-Ingersoll-Ross\ model} (CIR-model) introduced by Cox et al.\
(1985) for the term structure of interest rates. The above equation was
also studied in Ikeda and Watanabe (1989) and Revuz and Yor (1991). The
basic theory of general CBI-processes was developed in Li (2011). The
appealing properties of this process are as follows:
 \begin{enumerate}

\item[(1)] The process stays nonnegative.

\item[(2)] It converges to a steady-state law with mean $a/b$,
    the so-called long-term value, with speed of adjustment $b$.

\item[(3)] The incremental variance is proportional to its current
    value.

 \end{enumerate}
However, the one-dimensional CIR-model doesn't describe the connection
among interest rates of different currencies. In order to give more
objective description to the financial environment, we need to deal with
multifactor CIR-model or the multitype CBI-processes. In order to make the
presentation easier, we just analyze the two-type CBI-processes defined by
the following equation:
 \beqlb\label{1.2} \left\{\begin{array}{ccl}
dX_{1}(t)\ar=\ar(a_1-b_{11}X_1(t)
+b_{12}X_2(t))dt+\sigma_1\sqrt{X_1(t)}dB_1(t),\\
dX_{2}(t)\ar=\ar(a_2+b_{21}X_1(t)
-b_{22}X_2(t))dt+\sigma_2\sqrt{X_2(t)}dB_2(t),
\end{array}
\right.
 \eeqlb
where
$\theta=(a_1,a_2,b_{11},b_{12},b_{21},b_{22},\sigma_{1},\sigma_{2})\in
(0,+\infty)^3\times[0,\infty)^2\times(0,\infty)^3,$ and $B_i(t), i=1,2$ are independent standard Brownian
motions. The existence and uniqueness of the solution to (\ref{1.2}) were
proved in Ma (2012). We can rewrite (\ref{1.2}) into the vector form:
 \beqlb \label{1.3}
 dX_t=(A-BX_t)dt+\Sigma\sqrt{X_t}dW_t,
 \eeqlb
where
 \beqnn
 \begin{array}{rcl}
 X_t=\left(\begin{array}{c}X_1(t)\\
X_2(t)\end{array}\right), &A=\left(\begin{array}{c}a_1\\
a_2\end{array}\right),&
B=\left(\begin{array}{cc}b_{11}&-b_{12}\\
-b_{21}&b_{22}\end{array}\right),\\
\noalign{\vskip 3mm}
\Sigma=\left(\begin{array}{cc}\sigma_1&0\\
0&\sigma_2\end{array}\right),&\sqrt{X_t}=\left(\begin{array}{cc}\sqrt{X_1(t)}&0\\
0&\sqrt{X_2(t)}\end{array}\right),&W_t=\left(\begin{array}{c}B_1(t)\\
B_2(t)\end{array}\right).
\end{array}
 \eeqnn

We can use a special form of (\ref{1.2}) to describe the relations among
the interest rates of different currencies. In currency market, we assume
$X_1(t)$ is the interest rate of a very strong and influential currency,
and $X_2(t)$ is the interest rate of a less influential currency. The
situation can be described by the following stochastic equation:
 \beqnn
\left\{\begin{array}{ccl}
dX_{1}(t)\ar=\ar b_{11}(\frac{a_1}{b_{11}}-X_1(t))dt+\sigma_1\sqrt{X_1(t)}dB_1(t),\\
dX_{2}(t)\ar=\ar
b_{22}(\frac{a_2}{b_{22}}+\frac{b_{21}}{b_{22}}X_1(t)
-X_2(t))dt+\sigma_2\sqrt{X_2(t)}dB_2(t).
\end{array}
\right.
 \eeqnn
Here, the first equation gives the evolution of $X_1(t)$, which is
just a one-dimensional CIR- model. The second equation describes the
evolution of $X_2(t)$, which is affected not only by the random
noise, but also by $X_1(t)$. Specifically, the second coordination
$X_2(t)$ has the following properties:
 \begin{enumerate}

\item[(1)] It stays nonnegative.

\item[(2)] It converges to a steady-state law with mean
    $a_2/b_{22}+(b_{21}/b_{22}) (a_1/b_{11})$ (this can be easily got from (\ref{3.1}) with $t\rightarrow\infty$), the so-called long-term
    value, where the second term is contributed by $X_1(t)$.

\item[(3)] Its incremental variance is proportional to its current
    value.

 \end{enumerate}
If $b_{12}$ and $b_{21}$ are neither zero, then (\ref{1.2}) can account
for the fluctuations of the rate of two currencies that affect with each
other. Furthermore, (\ref{1.2}) can also be used to analyze the corporate
profitability and the market yield. For example,
 \begin{enumerate}

\item[(1)] In a perfectly competitive market, we can use (\ref{1.2})
    with $b_{12}=0$ to analyze the relationship between the corporate
    profitability and the yield of its corresponding industry.

\item[(2)] In an oligopoly market, we can use (\ref{1.2}) to describe
    the relationship among profitability of the different enterprises.

\item[(3)] However, in pure monopoly market, the one-dimensional
    CIR-model is suitable enough to analyze the corporate
    profitability or the market yield.
 \end{enumerate}

However, before using (\ref{1.3}) to solve the practical problems, we
need to estimate the parameters in the equation based on the historical
information. For the single-type CBI-processes, the approaches to
parameter estimation can be found in  Long-staff and Schwartz (1992) and
Bibby and S{\o}rensen (1995). Overbeck and Ryd\'{e}n (1997) also gave the
\textit{conditional least squares estimators} (CLSEs). Estimators of the
matrix of offspring means and the vector of stationary means in a
multitype GWI-process had been given in Quine and Durham (1977). For
multitype GWI-process, the \textit{weighted conditional least squares
estimator} (WCLSE) of the mean matrix was developed in Shete and Sriram
(2003). The asymptotically properties of CLSEs of GWI-processes with
general offspring laws were studied in Venkataraman (1982) and Wei and
Winnicki (1989). The asymptotics of CLSEs and WCLSEs of a stable
CIR-model was studied in Li and Ma (2013). It is well-known that the
CBI-processes are special examples of the affine Markov processes studied
in Duffie et al. (2003). The ergodicity and estimation of some different
two-dimensional affine processes were studied in Barczy et al. (2013a,
2013b, 2013c).

In this work, we give the CLSEs and the WCLSEs of the parameters in
(\ref{1.3}) using low frequency observations at equidistant time
points $\{k\Delta :k=0,1,\dots,n\}$ of a single realization
$\{X_t:t\geq 0\}$, where $X_t=(X_{1}(t),X_{2}(t))^{\mathrm{T}}$. For
simplicity, we take $\Delta=1$, but all the results presented below
can be modified to the general case. This is based on the
minimization of a sum of squared deviation about conditional
expectations developed in Klimko and Nelson (1978), who applied
their results to the CLSEs of the offspring and immigration means of
subcritical GWI-processes. Then, as Overbeck and Ryd\'{e}n (1997),
we shall study the consistency and the central limit theorems of
CLSEs and WCLSEs.

The paper is organized as follows. In Section~2, we give the
ergodicity of the two-type CBI-process, which is essentially
necessary for the study of the estimators. Section~3 is devoted to
the study of the CLSEs and WCLSEs of $(A,B)$ and $\Sigma$. The
consistency and asymptotic normality of the CLSEs and WCLSEs are
given in the Section~4. All the proofs are given in Section~5.


\section{ Multitype CBI-processes and ergodicity}
\setcounter{equation}{0}

 In this section, firstly we give the definition and a few properties of two-type
CBI-processes. In particular, we show that the solution of
(\ref{1.3}) is a two-type CBI-process. Secondly, we show that the
two-type CBI-process is ergodic under a weak condition. These
results are very important to study the consistency and asymptotic
normality of the estimators.

Let $\mbf{D}=[0, \infty)^2$. In view of (\ref{1.3}), we consider the
branching mechanisms $\phi_i, i=1,2$, with representation:
 \beqlb\label{2.1}
 \left\{\begin{array}{ccl}
 \phi_1(\lambda)\ar=\ar b_{11}\lambda_1-b_{12}\lambda_2+ \frac{\sigma^2_1}{2}
 \lambda_1^2,\cr
 \phi_2(\lambda)\ar=\ar -b_{21}\lambda_1+b_{22}\lambda_2+
 \frac{\sigma_2^2}{2}\lambda_2^2,
 \end{array}\right.
 \eeqlb
 where $\lambda=(\lambda_1,\lambda_2)^{\mathrm{T}}\in \mbf{D}$.
Next we give the definition of the two-type CBI-processes.
\begin{definition}
A Markov process $X_t=(X_1(t),X_2(t))$ in $\mbf{D}$ is called a two-type
continuous-state branching process with immigration (CBI-process), if it
has transition semigroup $(Q_t)_{t\geq 0}$ given by
 \beqlb \label{2.2}
\int_{\mbf{D}}e^{-\langle y,\lambda\rangle}Q_t(x,dy)=\exp\left\{-\langle
x,v_t(\lambda)\rangle-\int_0^t\langle A,v_s(\lambda)\rangle
ds\right\},\qquad \lambda,x\in \mbf{D},
 \eeqlb
where $A$ is given in (\ref{1.3}) and
$v_t(\lambda)=(v_1(t,\lambda),v_2(t,\lambda))^{\mathrm{T}}$ is the
unique solution of
 \beqlb \label{2.3}
 \left\{\begin{array}{ccl}
 \frac{\partial }{\partial t}v_1(t,\lambda)\ar=\ar-\phi_1(v_t(\lambda))= -b_{11}v_1(t,\lambda)+b_{12}v_2(t,\lambda)- \frac{\sigma^2_1}{2}v_1(t,\lambda)^2,\cr
 \frac{\partial }{\partial t}v_2(t,\lambda)\ar=\ar-\phi_2(v_t(\lambda))= b_{21}v_1(t,\lambda)-b_{22}v_2(t,\lambda)- \frac{\sigma^2_2}{2}v_2(t,\lambda)^2,
 \end{array}
 \right.
 \eeqlb
with the initial condition $v_0(\lambda)=\lambda \in \mbf{D}.$
\end{definition}

By Theorem~2.3 in Ma (2012), there is a unique non-negative weak
solution to (\ref{1.3}) and the solution is a CBI-process with
transition semigroup $(Q_t)_{t\geq0}$ defined by (\ref{2.2}). He
also showed that there is a unique non-negative strong solution to
(\ref{1.3}).

With the conclusions above, The following theorem gives a necessary
and sufficient condition for the ergodicity of the semigroup
$(Q_t)_{t\geq0}$.

Let $\kappa=\frac{b_{12}b_{21}}{b_{11}b_{22}}$. Then we get the following
conclusion, the proof will be given in Section~5.

 \begin{theorem}
Suppose both of the eigenvalues of $B=\bigg(\begin{array}{cc}b_{11}
&-b_{12}\cr -b_{21} &b_{22}\end{array}\bigg)$ are positive, i.e.
$\kappa<1$. Then the transition semigroup $(Q_t)_{t\geq0}$ has the
unique stationary distribution $Q_\infty$ given by
 \beqlb \label{2.15}
\int_{\mbf{D}}e^{-\langle
y,\lambda\rangle}Q_\infty(dy)=\exp\left\{-\int_0^\infty\langle
A,v_s(\lambda)\rangle ds \right\}.
 \eeqlb
Moreover, for every $x\in \mbf{D}$, we have $Q_t(x,\cdot)\rightarrow
Q_\infty$ by weak convergence as $t\rightarrow \infty.$
 \end{theorem}

With Theorem~2.2 the following corollary can be easily proved like
the proof of Theorem~2.7 in Li and Ma (2013).
 \begin{corollary}
Suppose the conditions of Theorem~2.2 are satisfied, then $X_t$ is
mixing and it's tail $\sigma$-algebra is trivial, e.g., Durrett
(2010).
 \end{corollary}

The results of this paper on the asymptotics of the estimators will
be derived under the assumption that the eigenvalues of $B$ are all
positive and $X_0$ is distributed according to the stationary law.
By Birkhoff's ergodic theorem, $X_t$ is a stationary and ergodic
process, but by a fairly simple (continuous time) coupling argument
it can be seen that they are valid for arbitrary initial
distributions.

\section{CLSEs and WCLSEs of $(A,B)$ and $\Sigma$}
\setcounter{equation}{0}

Recall $\Delta=1$. Let $\{\mathscr{F}_k=\sigma(\{X_0,\dots,X_k\}):k=0, 1,
2, \cdots\}$. In this section we give the CLSEs and WCLSEs of the drift
coefficients $(A,B)$ and the diffusion coefficients $\Sigma$ based on the
observations $\{X_k:k=0,\dots,n\}$.

Following Klimko and Nelson (1978) and Overbeck and Ryd\'{e}n
(1997), we first define the CLSE of $(A,B)$. The basic ideas are
explained as follows. By applying It\^{o}'s formula to (\ref{1.3}),
for any $t\ge r\ge 0$ we have
 \beqlb\label{3.1}
 X_t= e^{-B(t-r)}X_r + \int_r^t e^{-B(t-s)}A ds + \int_r^t e^{-B(t-s)}\Sigma \sqrt{X_s}dW_s.
 \eeqlb
Apparently, $\int_r^t e^{-B(t-s)}\Sigma \sqrt{X_s}dW_s$ is a local martingale with respect to $\{\mathscr{F}_t\}$. For
 \beqnn
 \mbf{E}\bigg[\int_0^te^{-B(t-s)}\left(\begin{array}{cc}\sigma_1^2X_1(t)&0\\0&\sigma_2^2X_2(t)\end{array}\right)e^{-B^{\mathrm{T}}(t-s)}dt\bigg]<\infty,\quad \forall\ t>0,
 \eeqnn
then  $\int_r^t e^{-B(t-s)}\Sigma \sqrt{X_s}dW_s$  is a martingale.

Let $\mbf{I}$ be the identical matrix and define
 \beqlb\label{3.2}
\rho= B^{-1}(\mathbf{I}-e^{-B})A, \quad \gamma= e^{-B}.
 \eeqlb
From (\ref{3.1}) and (\ref{3.2}) we can easily obtain the stochastic
regressive equation
 \beqlb \label{3.3}
 X_k= \rho+\gamma X_{k-1}+\varepsilon_k,
 \eeqlb
where
 $$\varepsilon_k=\int_{k-1}^{k} e^{-B(k-s)}\Sigma
\sqrt{X_{s}}dW_{s}.$$

Apparently, $\{\varepsilon_k:k\geq 0\}$ is a martingale differential
sequence with respect to $\{\mathscr{F}_k:k\geq 0\}$. Let
$m(x;A,B)=\rho+\gamma x$, then
 $$\mathbf{E}_\theta[X_k|X_{k-1}]=m(X_{k-1};A,B).$$

Suppose $g(x)$ is a Borel measurable function on $\mbb{R}^2$, satisfying
that
 $$\mathbf{E}(|g(X_0)|^2)<\infty,$$
let $g_k=g(X_{k-1})$, where $k=1,\dots,n$. The WCLSE of $(A,B)$ can
be given by minimizing the sum of squares
 \beqlb\label{3.4}
\sum_{k=1}^{n} g_k\varepsilon_k^{\mathrm{T}}\varepsilon_k
 =
\sum_{k=1}^{n} g_k(X_k-\rho-\gamma
X_{k-1})^{\mathrm{T}}(X_k-\rho-\gamma X_{k-1}).
 \eeqlb
 In particular, the estimators of $(\rho, \gamma)$ are given by
 \beqlb\label{3.5}
\hat{\rho}_n=\bar{g}_n\left(\bar{X}_n-\hat{\gamma}_n\tilde{X}_n\right)
 \eeqlb
 and
 \beqlb\label{3.6}
\hat{\gamma}_n=(\bar{T}_{1,n}-\tilde{T}_{1,n})\cdot(\bar{T}_{2,n}-\tilde{T}_{2,n})^{-1},
 \eeqlb
 where
 \begin{align*}
 \bar{X}_n&=\frac{1}{n}\sum_{k=1}^{n}g_kX_k,\quad \tilde{X}_n=\frac{1}{n}\sum_{k=1}^{n}g_kX_{k-1},\quad
\bar{g}_n=\frac{1}{n}\sum_{k=1}^{n}g_k, \\
\bar{T}_{1,n}&=\frac{1}{n}\sum_{k=1}^{n}(g_kX_k-\bar{X}_n)(g_kX_{k-1}-\tilde{X}_n)^{\mathrm{T}},\cr
\tilde{T}_{1,n}&=\frac{1}{n}\sum_{k=1}^{n}(g_k-\bar{g}_n)(g_kX_kX_{k-1}^{\mathrm{T}}-\frac{1}{n}\sum_{k=1}^{n}g_kX_kX_{k-1}^{\mathrm{T}}),\cr
\bar{T}_{2,n}&=\frac{1}{n}\sum_{k=1}^{n}(g_kX_{k-1}-\tilde{X}_n)(g_kX_{k-1}-\tilde{X}_n)^{\mathrm{T}},\\
\tilde{T}_{2,n}&=\frac{1}{n}\sum_{k=1}^{n}(g_k-\bar{g}_n)(g_kX_{k-1}X_{k-1}^{\mathrm{T}}-\frac{1}{n}\sum_{k=1}^{n}g_kX_{k-1}X_{k-1}^{\mathrm{T}}).
 \end{align*}

 According to  (\ref{3.2}), we
can easily get the
CLSE of $(A,B)$ by simple calculation,
 \beqnn
 \hat{A}_n\ar=\ar (\mbf{I}-\hat{\gamma}_{n})^{-1}\hat{B}_n\hat{\rho}_{n}, \\
 \hat{B}_n\ar=\ar-\ln \hat{\gamma}_{n},
 \eeqnn
Where $\ln\hat{\gamma}_{n}=\sum\limits_{k=0}^{\infty}(-1)^{k-1}\frac{1}{k}(\hat{\gamma}_{n}-1)^k$. From the proof of Theorem~4.1, $\hat{\gamma}_{n}\overset{\rm a.s.}\longrightarrow\gamma$ and the spectral radius of $\gamma$ is less than one, so $\hat{B}_n$ is reasonable.

We now turn to the estimation of $\Sigma$. For the connection
between $\Sigma$ and $(\sigma_1, \sigma_2)$, we just need to study
$(\sigma_1,\sigma_2).$ For this, we need to compute
 \beqlb\label{3.7}
v(X_{k-1};\theta)=\mathbf{E}_{\theta}[(X_k-m(X_{k-1};A,B))(X_k-m(X_{k-1};A,B))^{\mathrm{T}}|X_{k-1}].
 \eeqlb
 By applying It\^{o}'s formula to $X_tX_t^{\mathrm{T}}$, for
any $t\ge r\ge 0$ we have
 \beqnn X_tX_t^{\mathrm{T}}\ar=\ar
X_rX_r^{\mathrm{T}}
+\int_r^t\big[X_sA^{\mathrm{T}}+AX_s^{\mathrm{T}}+\Sigma\sqrt{X_s}\sqrt{X_s}^{\mathrm{T}}\Sigma^{\mathrm{T}}-(X_sX_s^{\mathrm{T}}B^{\mathrm{T}}-BX_sX_s^{\mathrm{T}})\big]ds\cr
\ar\ar+ \int_r^t \Sigma \sqrt{X_s}dW_sX_s^{\mathrm{T}}+\int_r^tX_s
dW_s^{\mathrm{T}} \sqrt{X_s}^{\mathrm{T}}\Sigma^{\mathrm{T}},
 \eeqnn
 where
 \beqnn
\Sigma\sqrt{X_s}\sqrt{X_s}^{\mathrm{T}}\Sigma^{\mathrm{T}}=\left(
 \begin{array}{cc}\sigma_1^2X_1(s)&0\\
 0&\sigma_2^2X_2(s)\end{array}\right).
 \eeqnn
Let $f(s)=\mathbf{E}_{\theta}[X_{k-1+s}|X_{k-1}]$ and
$h(s)=\mathbf{E}_{\theta}[X_{k-1+s}X_{k-1+s}^{\mathrm{T}}|X_{k-1}],\,s\geq
0$. Then for any $t\geq 0$
 \beqlb\label{3.8}
 h(t)= X_{k-1}X_{k-1}^{\mathrm{T}}+\int_0^t\Big\{f(s)A^{\mathrm{T}}+Af(s)^{\mathrm{T}}-(h(s)B^{\mathrm{T}}-Bh(s))+\left(
 \begin{array}{cc}\sigma_1^2f_1(s)&0\\
 0&\sigma_2^2f_2(s)
 \end{array}
 \right)\Big\}ds.
 \eeqlb
We can get $h(t)$ by solving (\ref{3.8}),
 \beqnn
h(t)\ar=\ar
e^{-Bt}X_{k-1}X_{k-1}^{\mathrm{T}}e^{-B^{\mathrm{T}}t}+\int_0^t
e^{-B(t-s)}(f(s)A^{\mathrm{T}}+Af(s)^{\mathrm{T}})e^{-B^{\mathrm{T}}(t-s)}ds\cr
\ar\ar\quad+\int_0^te^{-B(t-s)}\left(
 \begin{array}{cc}\sigma_1^2f_1(s)&0\\
 0&\sigma_2^2f_2(s)
 \end{array}\right)e^{-B^{\mathrm{T}}(t-s)}ds.
 \eeqnn
For
$$f(t)=e^{-Bt}X_{k-1} + \int_0^t
e^{-B(t-s)}Ads,$$ we can easily get
 \beqnn
f(t)f(t)^{\mathrm{T}}=e^{-Bt}X_{k-1}X_{k-1}^{\mathrm{T}}e^{-B^{\mathrm{T}}t}+\int_0^t
e^{-B(t-s)}(f(s)A^{\mathrm{T}}+Af(s)^{\mathrm{T}})e^{-B^{\mathrm{T}}(t-s)}ds.
 \eeqnn
  Thus
 \beqlb \label{3.9}
  v(X_{k-1};\theta)= h(1)-f(1)f(1)^{\mathrm{T}} =\sigma_1^2\eta_1+\sigma_2^2\eta_2,
 \eeqlb
where
 \beqnn
 \eta_1\ar=\ar\int_0^{1}e^{-B(1-s)}\left(
 \begin{array}{cc}f_1(s)&0\\
 0&0\end{array}\right)e^{-B^{\mathrm{T}}(1-s)}ds,\\
 \eta_2\ar=\ar\int_0^{1}e^{-B(1-s)}\left(
 \begin{array}{cc}0&0\\
 0&f_2(s)\end{array}\right)e^{-B^{\mathrm{T}}(1-s)}ds.
 \eeqnn

Before giving the estimator, we give the following important definition.
 \begin{definition}
Assume $A=(a_1,\dots,a_m)$ is $n\times m$ matrix, where
$\{a_i:i=1,\dots,m\}$ are  $n\times 1$ vectors, and define a
operator {\rm Vec} from a space of matrixes to a space of vectors,
i.e.
 \beqnn
 \mathrm{Vec}(A)=\left(\begin{array}{c}a_1\\ \vdots \\
a_m\end{array}\right).
 \eeqnn
 \end{definition}

Now we consider (\ref{3.9}) as a regression equation with
$(X_k-m(X_{k-1};A,B))(X_k-m(X_{k-1};A,B))^{\mathrm{T}}$
 being the response, $\sigma_1^2\eta_1+\sigma_2^2\eta_2$ being the predictor, and $(\sigma_1^2,\sigma_2^2)$ being the unknown. Recall $\{g_k:k=1,\dots,n\}$. The WCLSE of $(\sigma_1^2,\sigma_2^2)$ can be the minimizer of the following sum of squares, here we must of course substitute $(A,B)$ by the corresponding estimates
 \beqnn
\sum_{k=1}^ng_k\mathrm{Vec}(Z_k-(\sigma_1^2\hat{\eta}_{1,n}+\sigma_2^2\hat{\eta}_{2,n}))^{\mathrm{T}}\mathrm{Vec}(Z_k-(\sigma_1^2\hat{\eta}_{1,n}+\sigma_2^2\hat{\eta}_{2,n})),
 \eeqnn
where $(\hat{\eta}_{1,n},\hat{\eta}_{2,n})$ is $(\eta_1,\eta_2)$ with
$(A,B)$ substituted by $(\hat{A}_n,\hat{B}_n)$ and
 \beqnn
Z_k=(X_k-m(X_{k-1};\hat{A}_n,\hat{B}_n))(X_k-m(X_{k-1};\hat{A}_n,\hat{B}_n))^{\mathrm{T}}.
 \eeqnn
In particular, those estimators of $\sigma_1^2$ and $\sigma_2^2$
are given by
 \beqlb\label{3.10}
\hat{\sigma}^{2}_{1,n}=\frac{\varphi_{11,n}-\varphi_{12,n}}{\psi_n}
 \eeqlb
and
 \beqlb\label{3.11}
\hat{\sigma}^{2}_{2,n}=\frac{\varphi_{21,n}-\varphi_{22,n}}{\psi_n},
 \eeqlb
where $(\hat{\eta}_{1,n},\hat{\eta}_{2,n})$ is $(\eta_1,\eta_2)$
with $(A,B)$ substituted by
$(\hat{A}_n,\hat{B}_n)$
 \beqlb\label{3.12}
\varphi_{11,n}\ar=\ar\frac{1}{n}\sum_{k=1}^{n}g_k\mathrm{Vec}(Z_k)^{\mathrm{T}}\mathrm{Vec}(\hat{\eta}_{1,n})\mathrm{Vec}(\hat{\eta}_{2,n})^{\mathrm{T}}\mathrm{Vec}(\hat{\eta}_{2,n}),\\
\varphi_{12,n}\ar=\ar\frac{1}{n}\sum_{k=1}^{n}g_k\mathrm{Vec}(Z_k)^{\mathrm{T}}\mathrm{Vec}(\hat{\eta}_{2,n})\mathrm{Vec}(\hat{\eta}_{1,n})^{\mathrm{T}}\mathrm{Vec}(\hat{\eta}_{2,n}),\\
\varphi_{21,n}\ar=\ar\frac{1}{n}\sum_{k=1}^{n}g_k\mathrm{Vec}(Z_k)^{\mathrm{T}}\mathrm{Vec}(\hat{\eta}_{2,n})\mathrm{Vec}(\hat{\eta}_{1,n})^{\mathrm{T}}\mathrm{Vec}(\hat{\eta}_{1,n}),\\
\varphi_{22,n}\ar=\ar\frac{1}{n}\sum_{k=1}^{n}g_k\mathrm{Vec}(Z_k)^{\mathrm{T}}\mathrm{Vec}(\hat{\eta}_{1,n})\mathrm{Vec}(\hat{\eta}_{1,n})^{\mathrm{T}}\mathrm{Vec}(\hat{\eta}_{2,n})
 \eeqlb
and
 \beqlb\label{3.16}
\psi_n=\bar{g}_n\left[\mathrm{Vec}(\hat{\eta}_{1,n})^{\mathrm{T}}\mathrm{Vec}(\hat{\eta}_{1,n})\mathrm{Vec}(\hat{\eta}_{2,n})^{\mathrm{T}}\mathrm{Vec}(\hat{\eta}_{2,n})-\left(\mathrm{Vec}(\hat{\eta}_{1,n})^{\mathrm{T}}\mathrm{Vec}(\hat{\eta}_{2,n})\right)^2\right].
 \eeqlb

If choose $g(x)\equiv1$, we will get the CLSEs respectively. So in the later of this paper, we just need to deal with WCLSEs. If the observations of $X_t$ are too large,  we can
choose $g(x)=\frac{1}{1+|x|}$, where $|x|$ is the norm of $x$, this WCLSEs are usually discussed in other papers for its good properties.


\section{Consistency and asymptotic normality of $(\hat{A}_n,\hat{B}_n)$ and $(\hat{\sigma}^{2}_{1,n},\hat{\sigma}^{2}_{2,n})$}
\setcounter{equation}{0}

In this section we devote to show that $(\hat{A}_n,\hat{B}_n)$ and
$(\hat{\sigma}^{2}_{1,n},\hat{\sigma}^{2}_{2,n})$ are strongly
consistent. Further, we also analyze the central limit theorem of
the WCLSEs. All the proofs will be given in Section~5.

 \begin{theorem}
The estimator $(\hat{A}_n,\hat{B}_n)$ is strongly consistent, i.e.
$(\hat{A}_n,\hat{B}_n)\overset{\rm a.s.}\longrightarrow(A,B)$ as
$n\rightarrow \infty$.
 \end{theorem}

 \begin{remark}
 The estimators $\hat{A}_n,$ and $\hat{B}_n$ derived solely from
the conditional mean function are robust against misspecification of the
diffusion term in (\ref{1.3}).
 \end{remark}

 \begin{remark}
In the estimation of $(A,B)$, the diffusion term in (\ref{1.3}) can be replaced by $\sigma(X_t)dW_t$, where $\sigma(\cdot)$ is an arbitrary function such that the induced stationary distribution of $\{X_t\}$ has finite second moment.
 \end{remark}

 \begin{theorem}
The estimator $(\hat{\sigma}^{2}_{1,n},\hat{\sigma}^{2}_{2,n})$ are
strongly consistent, i.e.
$(\hat{\sigma}^{2}_{1,n},\hat{\sigma}^{2}_{2,n})\overset{\rm
a.s.}\longrightarrow(\sigma_1^2,\sigma_2^2)$ as $n\rightarrow
\infty$.
 \end{theorem}

 \begin{remark}
From the proof of Theorem~4.4 in Section~5, we can see any weakly
consistent estimator of $(A,B)$ gives weakly consistent estimator of
$\Sigma$.
 \end{remark}

Let
$\hat{\theta}_n=(\hat{a}_{1,n},\hat{a}_{2,n},\hat{b}_{11,n},\hat{b}_{12,n},\hat{b}_{21,n},\hat{b}_{22,n},\hat{\sigma}_{1,n},\hat{\sigma}_{2,n})$
is the WCLSE of $\theta$ given in the last section. Now we analyze
its asymptotic normality.

According to the argument above, we know $\hat{\theta}_n$ is the
unique solution of the following equation:
 \beqlb\label{4.1}
G_n(\theta)=\sum_{k=1}^{n}g_k\big\{w_0(X_{k-1};\theta)(X_k-\rho-\gamma
X_{k-1})+w_1(X_{k-1};\theta)\mathrm{Vec}(Z_k-(\sigma_1^2\eta_1+\sigma^2_2\eta_2))\big\}=0,
 \eeqlb
where
 \beqnn
w_0(x;\theta)=\left(\begin{array}{cc}\frac{\partial\rho_1}{\partial
a_1}&\frac{\partial\rho_2}{\partial a_1}\\
\frac{\partial\rho_1}{\partial
a_2}&\frac{\partial\rho_2}{\partial a_2} \\
\frac{\partial\rho_1}{\partial b_{11}}+(\frac{\partial
\gamma_{11}}{\partial b_{11}},\frac{\partial \gamma_{12}}{\partial
b_{11}})x & \frac{\partial\rho_2}{\partial b_{11}}+(\frac{\partial
\gamma_{21}}{\partial b_{11}},\frac{\partial
\gamma_{22}}{\partial b_{11}})x\\
\frac{\partial\rho_1}{\partial b_{12}}+(\frac{\partial
\gamma_{11}}{\partial b_{12}},\frac{\partial \gamma_{12}}{\partial
b_{12}})x & \frac{\partial\rho_2}{\partial b_{12}}+(\frac{\partial
\gamma_{21}}{\partial b_{12}},\frac{\partial
\gamma_{22}}{\partial b_{12}})x\\
\frac{\partial\rho_1}{\partial b_{21}}+(\frac{\partial
\gamma_{11}}{\partial b_{21}},\frac{\partial \gamma_{12}}{\partial
b_{21}})x & \frac{\partial\rho_2}{\partial b_{21}}+(\frac{\partial
\gamma_{21}}{\partial b_{21}},\frac{\partial
\gamma_{22}}{\partial b_{21}})x\\
\frac{\partial\rho_1}{\partial b_{22}}+(\frac{\partial
\gamma_{11}}{\partial b_{22}},\frac{\partial \gamma_{12}}{\partial
b_{22}})x & \frac{\partial\rho_2}{\partial b_{22}}+(\frac{\partial
\gamma_{21}}{\partial b_{22}},\frac{\partial
\gamma_{22}}{\partial b_{22}})x\\
0&0\\
0&0
\end{array}\right)
 \eeqnn
and
 \beqnn w_1(x;\theta)=\left(\begin{array}{cccccccc} 0,&0,&0,&0,&0,&0,&
\mathrm{Vec}(\eta_1),&
\mathrm{Vec}(\eta_2)
\end{array}\right)^{\mathrm{T}}.
 \eeqnn

The key to the analysis of $\hat{\theta}_n$ is that $G_n(\theta)$ is
a $\mathbb{P}_{\theta}$-martingale with respect to
$\{\mathscr{F}_n:n\geq 0\}.$ It's obvious, so we will not give the
proof.

Let
 $$\mu_3(x;\theta)=\mathbf{E}_{\theta}[(X_k-m(X_{k-1};\theta))\mathrm{Vec}(Z_k-v(x;\theta))^{\mathrm{T}}|X_{k-1}=x]$$
and
 $$\mu_4(x;\theta)=\mathbf{E}_{\theta}[\mathrm{Vec}(Z_k-v(x;\theta))\mathrm{Vec}(Z_k-v(x;\theta))^{\mathrm{T}}|X_{k-1}=x].$$
\begin{theorem}
 We have $\sqrt{n}(\hat{\theta}_n-\theta)\overset{\rm
d}\longrightarrow\mathcal {N}(0,V^{-1}WV^{-T})$ as
$n\rightarrow\infty$,
 where
 \beqlb \label{4.2}
W\ar=\ar\mathrm{E}_{\theta}\Big[g_1^2\big\{w_0(X_{0};\theta)v(X_{0};\theta)w_0^{\mathrm{T}}(X_{0};\theta)+w_1(X_{0};\theta)\mu_4(X_{0};\theta)w_1^{\mathrm{T}}(X_{0};\theta)\cr
\noalign{\vskip 3mm}
\ar\ar+w_0(X_{0};\theta)\mu_3(X_{0};\theta)w_1^{\mathrm{T}}(X_{0};\theta)+w_1(X_{0};\theta)\mu_3^{\mathrm{T}}(X_{0};\theta)w_0^{\mathrm{T}}(X_{0};\theta)\big\}\Big]
 \eeqlb
and
 \beqlb\label{4.3}
V=-\mathrm{E}_{\theta}\Big[g_1\big\{w_0(X_{0};\theta)\big(\frac{\partial
m(X_{0};\theta)}{\partial
\theta}\big)+w_1(X_{0};\theta)\big(\frac{\partial }{\partial
\theta}\mathrm{Vec}(v(X_{0};\theta))\big)\big\}\Big].
 \eeqlb
\end{theorem}

By taking $g\equiv 1$ one can see that all the conclusions above also
hold for the CLSEs.

\section{Proofs}

\setcounter{equation}{0}

In this section, we will give the proofs for the theorems in
Section~2 and 4

\noindent\textbf{Proof of Theorem~2.2} Let $u(t,\lambda)$ be the
unique solution of the linear equation:
 \beqlb \label{5.1}
\left\{\begin{array}{l}\frac{\partial}{\partial
t}u_1(t,\lambda)=-b_{11}u_1(t,\lambda)+b_{12}u_2(t,\lambda),\\
\frac{\partial}{\partial
t}u_2(t,\lambda)=b_{21}u_1(t,\lambda)-b_{22}u_2(t,\lambda),\end{array}\right.
 \eeqlb
with the initial condition $u(0,\lambda)=\lambda\in \mbf{D}.$ By applying
the  comparison theorem to (\ref{2.3}) and (\ref{5.1}), we have
$v_i(t,\lambda)\leq u_i(t,\lambda)$ for $t\geq 0$, where $i=1,2$. By
solving (\ref{5.1}), we get $u(t,\lambda) = e^{-Bt}\lambda$. Let
$\xi_1,\xi_2>0$ be the eigenvalues of $B$. If $\xi_1=\xi_2$, then
$b_{11}=b_{22}>0$, $b_{12}=b_{21}=0$ and
 \beqlb\label{5.2}
u(t,\lambda)=\left(\begin{array}{c}e^{-b_{11}t}\lambda_1\cr
e^{-b_{22}t}\lambda_2\end{array}\right).
 \eeqlb If $\xi_1\neq\xi_2$,
let $L_i=\{\ell:(B-\xi_i\mbf{I})\ell=0\}$, where $i=1,2$. Then
$V_1,V_2$ are subspaces of $\mbb{R}^2$. For any $\lambda\in
\mbf{D}$, there exit $\ell_1\in L_1$ and $\ell_2\in L_2$, such that
$\lambda=\ell_1+\ell_2$ and
 \beqlb\label{5.3}
(B-\xi_i\mbf{I})^l\ell_i=0\quad l\geq1,\ i=1,2.
 \eeqlb
 By (\ref{5.3}), we have
 \beqlb\label{5.4}
e^{-Bt}\lambda=e^{-\xi_1t}\ell_1+e^{-\xi_2t}\ell_2.
 \eeqlb
By (\ref{5.2}) and (\ref{5.4}), we have
 \beqlb\label{5.5}
v_i(t,\lambda)\leq u_i(t,\lambda)\leq C_1(\lambda)e^{-(\xi_1\wedge\,
\xi_2)t} \quad \mbox{for some}\ C_1(\lambda)>0.
 \eeqlb
 Thus $\lim\limits_{t\rightarrow +\infty}v_i(t,\lambda)=0$. So for $\forall\,
x\in\mbf{D}$,
 \beqnn
\lim_{t\rightarrow +\infty}\int_{\mbf{D}}e^{-\langle
y,\lambda\rangle}Q_t(x,dy)=\lim_{t\rightarrow +\infty}e^{-\langle
x,v_t(\lambda)\rangle-\int_0^t\langle A,v_s(\lambda)\rangle
ds}=e^{-\int_0^{\infty}\langle A,v_s(\lambda)\rangle ds}.
 \eeqnn
By (\ref{5.5}), we have $\langle
A,v_s(\lambda)\rangle=a_1v_1(s,\lambda)+a_2v_2(s,\lambda)\leq
C_2(\lambda)e^{-(\xi_1\wedge\,\xi_2)t}$, for some $C_2(\lambda)>0$.
Then
 \beqnn \int_0^{\infty}\langle A,v_s(\lambda)\rangle
ds<\infty\quad \mbox{for}\ \forall~\lambda\in\mbf{D}.
 \eeqnn
 Thus, there is a unique distribution $Q_{\infty}(\cdot)$ satisfying
 \beqnn
\int_{\mbf{D}}e^{-\langle
y,\lambda\rangle}Q_\infty(dy)=\exp\left\{-\int_0^\infty\langle
A,v_s(\lambda)\rangle ds \right\}.
 \eeqnn\qed

%

\smallskip \noindent\textbf{Proof of Theorem~4.1} By the
relationship of $(\hat{A}_n,\hat{B}_n)$ and
$(\hat{\rho}_n,\hat{\gamma}_n)$, we just need to prove the
consistency of $(\hat{\rho}_n,\hat{\gamma}_n)$, i.e.
 $$(\hat{\rho}_n,\hat{\gamma}_n)\overset{\rm
a.s.}\longrightarrow(\rho,\gamma) \ \mbox{as}\ n\rightarrow\infty.$$
By ergodicity, we have
 \beqnn
\bar{T}_{1,n}\ar\overset{\rm
a.s.}\longrightarrow\ar\mathrm{Cov}(g_1X_1,g_1X_{0}) \ \mbox{as}\
n\rightarrow\infty
 \eeqnn
 and
 \beqnn \tilde{T}_{1,n}\ar\overset{\rm
a.s.}\longrightarrow\ar\mathrm{Cov}(g_1,g_1X_1X_{0}^{\mathrm{T}})\
\mbox{as}\ n\rightarrow\infty.
 \eeqnn
 For
  \beqnn
\mathrm{Cov}(g_1X_1,g_1X_{0})\ar=\ar\mathbf{E}_\theta[g_1^2\mathbf{E}_\theta(X_1|\mathscr{F}_{0})X_{0}^{\mathrm{T}}]-\mathbf{E}_\theta[g_1\mathbf{E}_\theta(X_1|\mathscr{F}_{0})]\mathbf{E}_\theta[g_1X_{0}^{\mathrm{T}}],\\
\mathrm{Cov}(g_1,g_1X_1X_{0}^{\mathrm{T}})\ar=\ar\mathbf{E}_\theta[g_1^2\mathbf{E}_\theta(X_1|\mathscr{F}_{0})X_{0}^{\mathrm{T}}]-\mathbf{E}_\theta[g_1]\mathbf{E}_\theta[g_1\mathbf{E}_\theta(X_1|\mathscr{F}_{0})X_{0}^{\mathrm{T}}],
  \eeqnn
 then we get
 \beqlb\label{5.6}
\bar{T}_{1,n}-\tilde{T}_{1,n}\overset{\rm
a.s.}\longrightarrow\gamma\{\mathbf{E}_\theta[g_1]\mathbf{E}_\theta[g_1X_{0}X_{0}^\mathrm{T}]-\mathbf{E}_\theta[g_1X_{0}]\mathbf{E}_\theta[g_1X_{0}^\mathrm{T}]\}\
\mbox{as}\ n\rightarrow\infty.
 \eeqlb
 Similarly, we can easily get
 \beqlb\label{5.7}
\bar{T}_{2,n}-\tilde{T}_{2,n}\overset{\rm
a.s.}\longrightarrow\mathbf{E}_\theta[g_1]\mathbf{E}_\theta[g_1X_{0}X_{0}^\mathrm{T}]-\mathbf{E}_\theta[g_1X_{0}]\mathbf{E}_\theta[g_1X_{0}^\mathrm{T}]\
\mbox{as}\ n\rightarrow\infty.
 \eeqlb
  By (\ref{3.6}), (\ref{5.6}) and (\ref{5.7}),
 $$\hat{\gamma}_n\overset{\rm
a.s.}\longrightarrow\gamma\ \mbox{as}\ n\rightarrow\infty.$$
Similarly, by ergodicity and (\ref{3.5}),
 $$\hat{\rho}_n\overset{\rm
a.s.}\longrightarrow\rho \ \mbox{as}\ n\rightarrow\infty.$$ \qed

\smallskip
\noindent\textbf{Proof of Theorem~4.4} By ergodicity and theorem~4.1, we can get
$$(\mathrm{Vec}(\hat{\eta}_{1,n}),\mathrm{Vec}(\hat{\eta}_{2,n}))\overset{\rm
a.s.}\longrightarrow(\mathrm{Vec}(\eta_1),\mathrm{Vec}(\eta_2))\
\mbox{as}\ n\rightarrow\infty.$$ Next we will prove $\psi_n$ and
$\varphi_{ij,n}$ converge almost surly, where $i,j=1,2$. Fix
$\theta.$ For $\forall \vartheta\in[0,\infty)^8$, let
 \beqnn h(x,y;\vartheta)\ar=\ar
g(x)\mathrm{Vec}((y-m(x;\vartheta))(y-m(x;\vartheta))^{\mathrm{T}})^{\mathrm{T}}
\mathrm{Vec}(\eta_1(\vartheta))\mathrm{Vec}(\eta_2(\vartheta))^{\mathrm{T}}\mathrm{Vec}(\eta_2(\vartheta))
 \eeqnn
 and $U\subset [0,+\infty)^8$ be a neighborhood of $\theta$
such that
 \beqnn\mathbf{E}_\theta[\sup_{\vartheta\in
U}|h(X_0,X_1;\vartheta)|]<\infty.
 \eeqnn
  By ergodicity, (\ref{3.12}) and Theorem~4.1
 \beqnn \varlimsup_{n\rightarrow
\infty}\varphi_{11,n} \leq\varlimsup_{n\rightarrow
\infty}\frac{1}{n}\sum_{k=1}^{n}\sup_{\vartheta\in
U}h(X_{k-1},X_k;\vartheta)=\mathbf{E}_\theta[\sup_{\vartheta\in
U}h(X_0,X_1;\vartheta)].
 \eeqnn
Let $U\downarrow \{\theta\}$. Then
 \beqlb\label{5.8}
\varlimsup_{n\rightarrow \infty}\varphi_{11,n}\leq
\mathbf{E}_\theta[h(X_0,X_1;\theta)].
 \eeqlb
 Similarly, we can prove
 \beqlb\label{5.9}
\varliminf_{n\rightarrow \infty}\varphi_{11,n}\geq
\mathbf{E}_\theta[h(X_0,X_1;\theta)].
 \eeqlb
By (\ref{5.8}) and (\ref{5.9}),
 \beqlb\label{5.10}
\varphi_{11,n}\ar\overset{\rm
a.s.}\longrightarrow\ar\sigma_1^2\mathbf{E}_\theta\big[g_1\mathrm{Vec}(\eta_1)^{\mathrm{T}}
\mathrm{Vec}(\eta_1)\mathrm{Vec}(\eta_2)^{\mathrm{T}}\mathrm{Vec}(\eta_2)\big]\cr\noalign{\vskip
3mm}
\ar\ar+\sigma_2^2\mathbf{E}_\theta[g_1\mathrm{Vec}(\eta_2)^{\mathrm{T}}
\mathrm{Vec}(\eta_1)\mathrm{Vec}(\eta_2)^{\mathrm{T}}\mathrm{Vec}(\eta_2)]\
\mbox{as}\ n\rightarrow\infty .
 \eeqlb Similarly, when
$n\rightarrow\infty$ we can get \
 \beqlb\label{5.11}
\varphi_{12,n}\overset{\rm
a.s.}\longrightarrow\sigma_1^2\mathbf{E}_\theta\big[g_1(\mathrm{Vec}(\eta_1)^{\mathrm{T}}
\mathrm{Vec}(\eta_2))^2\big]+\sigma_2^2\mathbf{E}_\theta\big[g_1\mathrm{Vec}(\eta_2)^{\mathrm{T}}
\mathrm{Vec}(\eta_1)\mathrm{Vec}(\eta_2)^{\mathrm{T}}\mathrm{Vec}(\eta_2)\big]
 \eeqlb
and
 \beqlb\label{5.12}
 \psi_n\overset{\rm
a.s.}\longrightarrow\mathbf{E}_\theta[g_1\mathrm{Vec}(\eta_1)^{\mathrm{T}}\mathrm{Vec}(\eta_1)\mathrm{Vec}(\eta_2)^{\mathrm{T}}\mathrm{Vec}(\eta_2)-g_1(\mathrm{Vec}(\eta_1)^{\mathrm{T}}\mathrm{Vec}(\eta_2))^2].
 \eeqlb
By (\ref{3.10}), (\ref{5.10}), (\ref{5.11}) and (\ref{5.12}),
 $$
\hat{\sigma}^{2}_{1,n}\overset{\rm a.s.}\longrightarrow\sigma_1^2 \
\mbox{as}\ n\rightarrow\infty.
 $$
Similarly, we can prove
 $$\hat{\sigma}^{2}_{2,n}\overset{\rm
a.s.}\longrightarrow \sigma_2^2 \ \mbox{as}\ n\rightarrow\infty.$$
\qed

\smallskip \noindent\textbf{Proof of Theorem~4.6}  Fix $\theta$. By
Theorem~4.1 and 4.4, making a Taylor expansion of $G_n$ about
$\theta$, it is enough to show that
 \beqlb \label{5.13}
G_n(\theta^{\,'})=G_n(\theta)+G^{\,'}_{n}(\theta)(\theta^{\,'}-\theta)+\frac{1}{2}G^{\,''}_{n}(\xi)(\theta^{\,'}-\theta)(\theta^{\,'}-\theta)^{\mathrm{T}},
 \eeqlb
where $\xi$ is between $\theta^{\,'}$ and $\theta$. Each element of
(\ref{5.13}) is
 \beqnn
 G_{n,i}(\theta^{\,'})\ar=\ar
G_{n,i}(\theta)+\sum_{l=1}^{8} \frac{\partial}{\partial
\theta_l}G_{n,i}(\theta)(\theta^{\,'}_l-\theta_l)+\frac{1}{2}\sum_{l=1}^{8}\sum_{j=1}^{8}\frac{\partial
^2}{\partial \theta_l\partial
\theta_j}G_{n,i}(\xi)(\theta^{\,'}_j-\theta_j)(\theta^{\,'}_l-\theta_l),
 \eeqnn
where\, $i=1,\dots,8$.

If choose $\theta^{\,'}=\hat{\theta}_n$ in the above equation, we
will get
 \beqlb\label{5.14}
-\frac{1}{\sqrt{n}}G_{n,i}(\theta)=\sum_{l=1}^{8}\sqrt{n}(\hat{\theta}_{n,l}-\theta_l)
\Big\{\frac{1}{n}\frac{\partial}{\partial
\theta_l}G_{n,i}(\theta)+\frac{1}{2n}\sum_{j=1}^{8}\frac{\partial
^2}{\partial \theta_l\partial
\theta_j}G_{n,i}(\xi)(\hat{\theta}_{n,j}-\theta_j)\Big\}.
 \eeqlb
 Let
 $$T_{n,l}=\sqrt{n}(\hat{\theta}_{n,l}-\theta_l),\quad
Y_{n,i}=-\frac{1}{\sqrt{n}}G_{n,i}(\theta),$$
$$ D_{n,il}=\frac{1}{n}\frac{\partial}{\partial
\theta_l}G_{n,i}(\theta)+\frac{1}{2n}\sum_{j=1}^{8}\frac{\partial
^2}{\partial \theta_l\partial
\theta_j}G_{n,i}(\xi)(\hat{\theta}_{n,j}-\theta_j).$$
 Then (\ref{5.14})
turns to be
 \beqlb\label{5.15}
  D_nT_n=Y_n,
 \eeqlb where
 $D_n=(D_{n,il})_{8\times8}$,
 $T_n=(T_{n,1},\dots,T_{n,8})^\mathrm{T}$,
 $Y_n=(Y_{n,1},\dots,Y_{n,8})^\mathrm{T}$.
Next we will prove the following results respectively,
 \beqlb
\label{5.16} D_n\overset{\rm a.s.}\longrightarrow V\ \mbox{as}\
n\rightarrow\infty
 \eeqlb
  and
 \beqlb \label{5.17}
 Y_n\ar\overset{\rm
d}\longrightarrow\ar \mathcal {N}(0,W)\ \mbox{as}\
n\rightarrow\infty.
 \eeqlb
By (\ref{5.15}), (\ref{5.16}) and (\ref{5.17}), we can get
 $$T_n\overset{\rm
d}\longrightarrow\mathcal {N}(0,V^{-1}WV^{-T})\ \mbox{as}\
n\rightarrow\infty.$$

{\bf(1)} Firstly, we prove (\ref{5.16}) holds.
Let
 \beqnn
h(x,y;\theta)=w_0(x;\theta)(y-m(x;\theta))+w_1(x;\theta)\mathrm{Vec}(z-v(x;\theta)),
 \eeqnn
where $z=(y-m(x;A,B))(y-m(x;A,B))^{\mathrm{T}}$. Then (\ref{4.1})
turns out to be
 \beqnn
G_n(\theta)=\sum_{k=1}^{n}g_kh(X_{k-1},X_k;\theta).
 \eeqnn
It's easy to prove that\ $\exists H(x,y)$ and $\theta\in U$, such
that
 \begin{enumerate}
\item[(i)]$|\frac{\partial^2}{\partial\theta_j\theta_l}h_i(x,y;\theta)|\leq
H(x,y)$ uniformly in $U$, for $i,j,l=1,\dots,8$.
\item[(ii)]$\mathrm{E}_\theta[|H(X_0,X_1)|]<\infty$.
 \end{enumerate}
Then
 \beqnn \left|\frac{1}{2n}\sum_{j=1}^{8}\frac{\partial
^2}{\partial \theta_l\partial
\theta_j}G_{n,i}(\xi)(\hat{\theta}_{n,j}-\theta_j)\right| \ar=\ar
\left|\frac{1}{2n}\sum_{j=1}^{8}(\hat{\theta}_{n,j}-\theta_j)\sum_{k=1}^{n}g_k\frac{\partial
^2}{\partial \theta_l\partial
\theta_j}h_i(X_{k-1};X_k;\xi)\right|\\
&\leq\ar\sum_{j=1}^{8}|\hat{\theta}_{n,j}-\theta_j|\frac{1}{2n}\sum_{k=1}^{n}|g_k|H(X_{k-1},X_k).
 \eeqnn
By Theorem~4.1, 4.4, and
 \beqnn
\frac{1}{2n}\sum_{k=1}^{n}|g_k|H(X_{k-1},X_k) &\overset{\rm
a.s.}\longrightarrow\frac{1}{2}\mathrm{E}_{\theta}[|g_1|H(X_{0},X_{m})]<\infty\
\mbox{as}\ n\rightarrow\infty,
 \eeqnn
we can get
 \beqnn \left|\frac{1}{2n}\sum_{j=1}^{8}\frac{\partial
^2}{\partial \theta_l\partial
\theta_j}G_{n,i}(\xi)(\hat{\theta}_{n,j}-\theta_j)\right|
&\overset{\rm a.s.}\longrightarrow 0\ \mbox{as}\ n\rightarrow\infty.
 \eeqnn
By (\ref{4.1}),
 \beqnn
\frac{1}{n}\frac{\partial}{\partial\theta_l}G_{n,i}(\theta)\ar=\ar
\frac{1}{n}\sum_{k=1}^{n}g_k\Big\{(\frac{\partial}{\partial\theta_l}w_{0,i}(X_{k-1};\theta))(X_k-m(X_{k-1};\theta))\cr
\ar\ar+w_{0,i}(X_{k-1};\theta)(-\frac{\partial}{\partial\theta_l}m(X_{k-1};\theta))\cr
\ar\ar+(\frac{\partial}{\partial\theta_l}w_{1,i}(X_{k-1};\theta))\mathrm{Vec}(Z_k-v(X_{k-1};\theta))^{\mathrm{T}}\cr
\ar\ar
+w_{1,i}(X_{k-1};\theta)\mathrm{Vec}\Big((-\frac{\partial}{\partial\theta_l}m(X_{k-1};\theta))(X_k-m(X_{k-1};\theta))^{\mathrm{T}}\cr
\ar\ar+(X_k-m(X_{k-1};\theta))(-\frac{\partial}{\partial\theta_l}m(X_{k-1};\theta))^{\mathrm{T}}-\frac{\partial}{\partial\theta_l}v(X_{k-1};\theta)\Big)\Big\},
 \eeqnn
where $\{w_{0,i},w_{1,i}:i=1,\dots,8\}$ are the row vectors of $w_0$
and $w_1$. By ergodicity, we calculate each part respectively, we
can get
 \beqnn\frac{1}{n}\frac{\partial}{\partial
\theta_l}G_{n,i}(\theta)\rightarrow V_{i,l},
 \eeqnn
where
 $$
V_{i,l}=-\mathrm{E}_{\theta}\Big[g_1\big[w_{0,i}(X_{0};\theta)\frac{\partial}{\partial\theta_l}m(X_{0};\theta)+w_{1,i}(X_{0};\theta)\mathrm{Vec}(\frac{\partial}{\partial\theta_l}v(X_{0};\theta))\big]\Big].
 $$
Thus we get (\ref{5.16}).

{\bf(2)} Secondly, we prove (\ref{5.17}) holds.
Recall $\{G_n(\theta)\}$ is a martingale, let
 \beqnn
\frac{1}{n}V_n\ar=\ar\frac{1}{n}\sum_{k=1}^{n}\mathrm{E}_{\theta}[(g_kh(X_{k-1},X_k;\theta))(g_kh(X_{k-1},X_k;\theta))^{\mathrm{T}}|\mathscr{F}_{k-1}]\cr
\ar=\ar\frac{1}{n}\sum_{k=1}^{n}g_k^2\big[w_0(X_{k-1};\theta)v(X_{k-1};\theta)w_0^{\mathrm{T}}(X_{k-1};\theta)
+w_1(X_{k-1};\theta)\mu_4(X_{k-1};\theta)w_1^{\mathrm{T}}(X_{k-1};\theta)\\
\ar\ar
+w_0(X_{k-1};\theta)\mu_3(X_{k-1};\theta)w_1^{\mathrm{T}}(X_{k-1};\theta)
+w_1(X_{k-1};\theta)\mu_3^{\mathrm{T}}(X_{k-1};\theta)w_0^{\mathrm{T}}(X_{k-1};\theta)\big].
 \eeqnn
By ergodicity, we can easily get
 $$
\frac{1}{n}V_n\rightarrow W\ \mbox{as}\ n\rightarrow \infty.
 $$
By stationarity, For $\forall \varepsilon >0$
 \beqnn
\lefteqn{\frac{1}{n}\sum_{k=1}^{n}\mathrm{E}_{\theta}\big[g_k^2h(X_{k-1},X_k;\theta)^{\mathrm{T}}h(X_{k-1},X_k;\theta)1_{\{g_k^2h(X_{k-1},X_k;\theta)^{\mathrm{T}}h(X_{k-1},X_k;\theta)>n\varepsilon^2\}}\big]}\cr
&&=\mathrm{E}_{\theta}[(g_1^2h(X_{0},X_1;\theta)^{\mathrm{T}}h(X_{0},X_1;\theta)1_{\{g_1^2h(X_{0},X_1;\theta)^{\mathrm{T}}h(X_{0},X_1;\theta)>n\varepsilon^2\}}]\rightarrow0\
\mbox{as}\ n\rightarrow\infty.
 \eeqnn
Thus, by the martingale central theorem (see, e.g., Durrett, 2010),
(\ref{5.17}) holds.

\qed


\end{document}